\documentclass[10pt,a4paper,twoside,leqno]{article}
\usepackage{amsfonts,amssymb}

\usepackage{amscd}

\vspace{0.5cm}

 4
\font\ebf=cmbx8
\font\erm=cmr8

\usepackage{graphicx}

\author{M.~Dziemia{\'n}czuk}
\title{On Cobweb posets tiling problem}

\newtheorem{defn}{Definition}
\newtheorem{theoremn}{Theorem}
\newtheorem{observen}{Observation}
\newtheorem{example}{Example}

\begin {document}


\begin{center}
\noindent {\bf\Large On Cobweb posets tiling problem}

\vspace{0.5cm}

\noindent {M. Dziemia\'nczuk}

\vspace{0.5cm}

\noindent {\erm Physics Department Bia\l ystok University (*)}

\noindent {\erm  Sosnowa 64, PL-15-887 Bia\l ystok, Poland}

\noindent {\erm e-mail: Maciek.Ciupa@gmail.com}

\noindent {\erm (*) former Warsaw University Division}
\end{center}

\vspace{1cm}

\noindent {\ebf SUMMARY}

\vspace{0.1cm}

\noindent {\small Kwa\'sniewski's cobweb posets uniquely represented by directed acyclic graphs are such a generalization of the Fibonacci tree that allows joint combinatorial interpretation for all of them under admissibility condition. This interpretation was derived in the source papers and it entailes natural enquieres already formulated therein.  In our note we response to one of those problems.  This is a tiling problem. Our observations on tiling problem include proofs of tiling's existence for some cobweb-admissible sequences. We show also that not all cobwebs admit tiling as defined below. }

\vspace{0.3cm}


\vspace{0.1cm}

\noindent AMS Classification Numbers: 06A07, 05C70, 05C75, 11B39.
\vspace{0.1cm}

\noindent  Presented at Gian-Carlo Polish Seminar:

\noindent \emph{http://ii.uwb.edu.pl/akk/sem/sem\_rota.htm}

\vspace{0.4cm}

\section{Introduction}

The source papers are \cite{1,2} from which indispensable definitions and notation are taken for granted as for example  (Kwa\'sniewski upside - down notation $n_F \equiv F_n$ being used for mnemonic reasons \cite{1,2,3}) :
$F-nomial$ coefficient:

\begin{displaymath}
	{n \choose k}_F = \frac{n_F \cdot (n-1)_F \cdot \ldots \cdot (n-k+1_F }
{ 1_F \cdot 2_F \cdot \ldots \cdot k_F } = \frac{ n_{F}^{\underline{k}} }{ k_F!}; \quad n_F \equiv F_n
\end{displaymath}

\begin{figure}[ht]
\begin{center}
	\includegraphics[width=100mm]{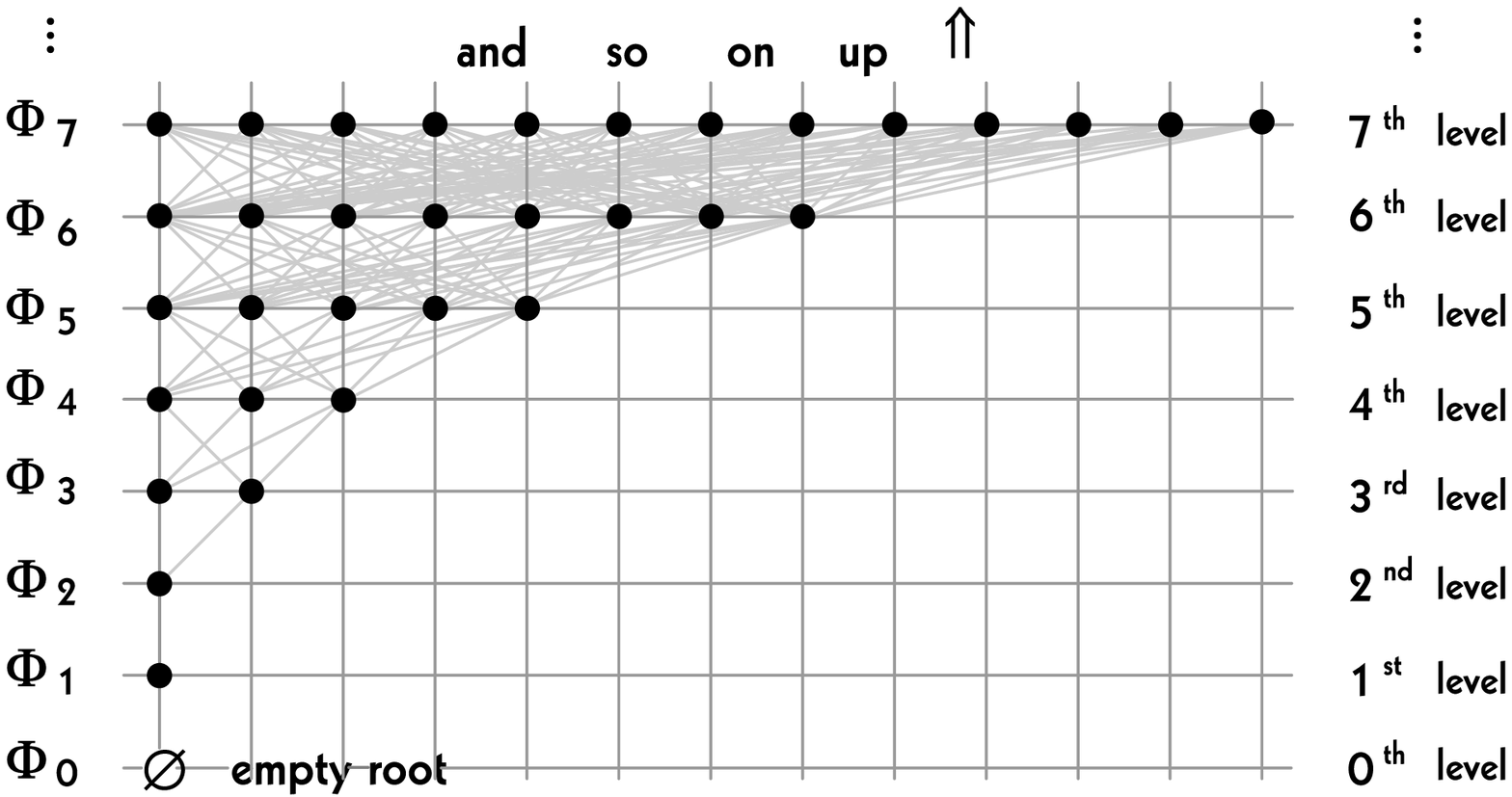}
	\caption{The $s$-th level in $N\times N_0$ \label{fig:akkfig0}}
\end{center}
\end{figure}

Nevertheless let us at first recall that cobweb poset in its original form \cite{1,2} is defined as a partially ordered graded infinite poset $\Pi = \langle P,\leq\rangle$, designated uniquely by any sequence of nonnegative integers $F = \{ n_F \}_{n \geq 0 }$ and it is represented as a directed acyclic graph (DAG) in the graphical display of its Hasse diagram. $P$ in $\langle P,\leq \rangle$ stays for   set of vertices while $\leq$ denotes partially ordered relation. See Figure \ref{fig:akkfig0}. and note  (quotation from \cite{2,1}):

\begin{quote}
One  refers to $\Phi_s$ as to  the set of vertices at the $s$-th level. The population of the $k$-th level ("\emph{generation}") counts  $k_F$ different member vertices for $k>0$ and one for $k=0$.
Here down (Fig. \ref{fig:akkfig0}) a disposal of vertices on $\Phi_k$  levels is visualized for the case of Fibonacci sequence. $F_0 = 0$ corresponds to the empty root $\{\emptyset\}$.
\end{quote}

In Kwa\'sniewski's cobweb posets' tiling problem one considers finite cobweb sub-posets for which we have finite number of layers $\langle \Phi_k \rightarrow \Phi_n \rangle$, where  $k\leq n$,  $k, n\in\mathbb{N}\cup \{0\} $ with exactly $k_j$  vertices on $\Phi_j$ level $k\leq j \leq n$. For $k=0$ the sub-posets $\langle \Phi_0 \rightarrow \Phi_n \rangle$ are named \emph{ prime cobweb posets} and these are those to be used - up to permutation of levels equivalence - as a block to partition finite cobweb sub-poset.

For the sake of combinatorial interpretation \cite{1,2}  a natural numbers valued sequence $F$ which determines a cobweb poset has to be  the so-called \emph{cobweb-admissible}.

\begin{figure}[ht]
\begin{center}
	\includegraphics[width=75mm]{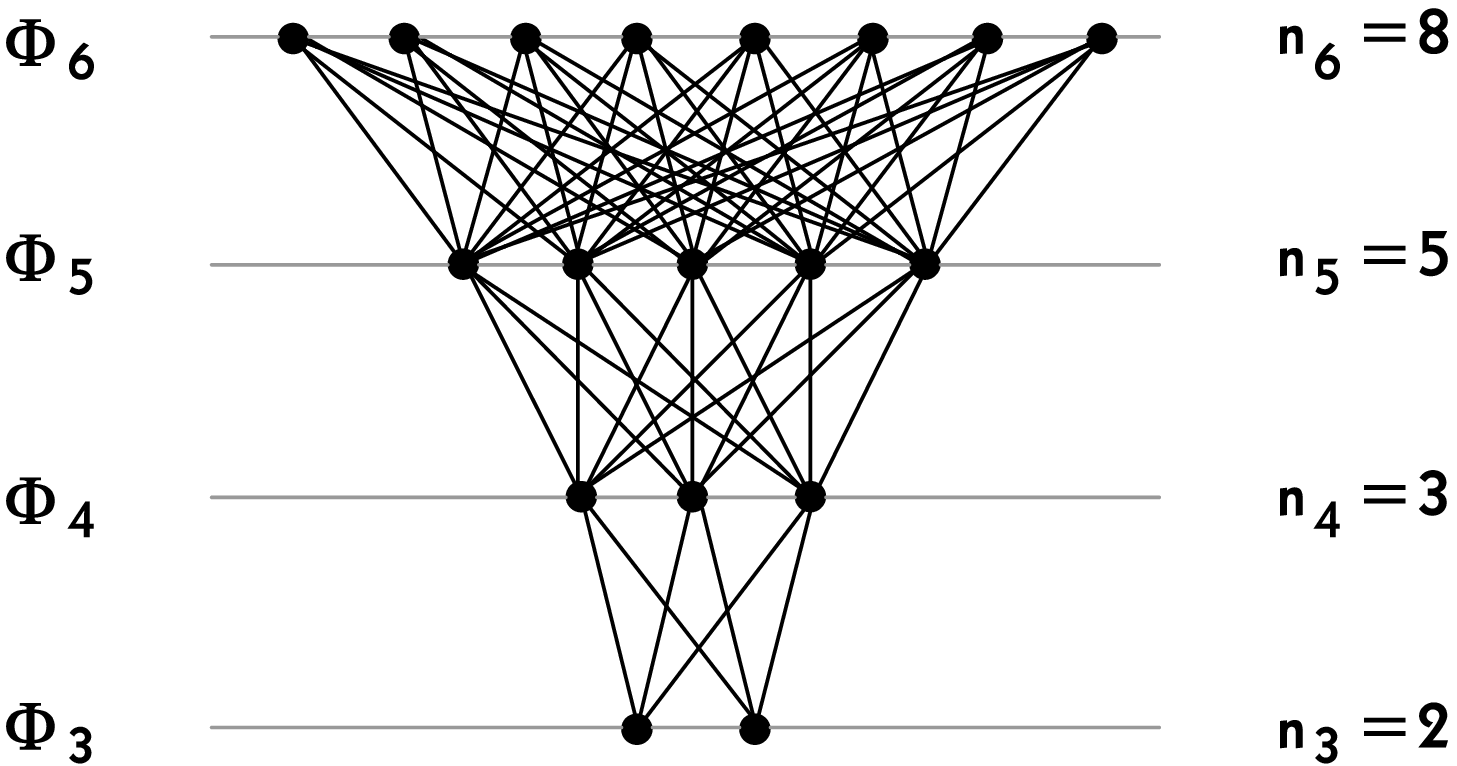}
	\caption{Display of four levels of Fibonacci numbers' finite Cobweb sub-poset}
\end{center}
\end{figure}

\begin{figure}[ht]
\begin{center}
	\includegraphics[width=75mm]{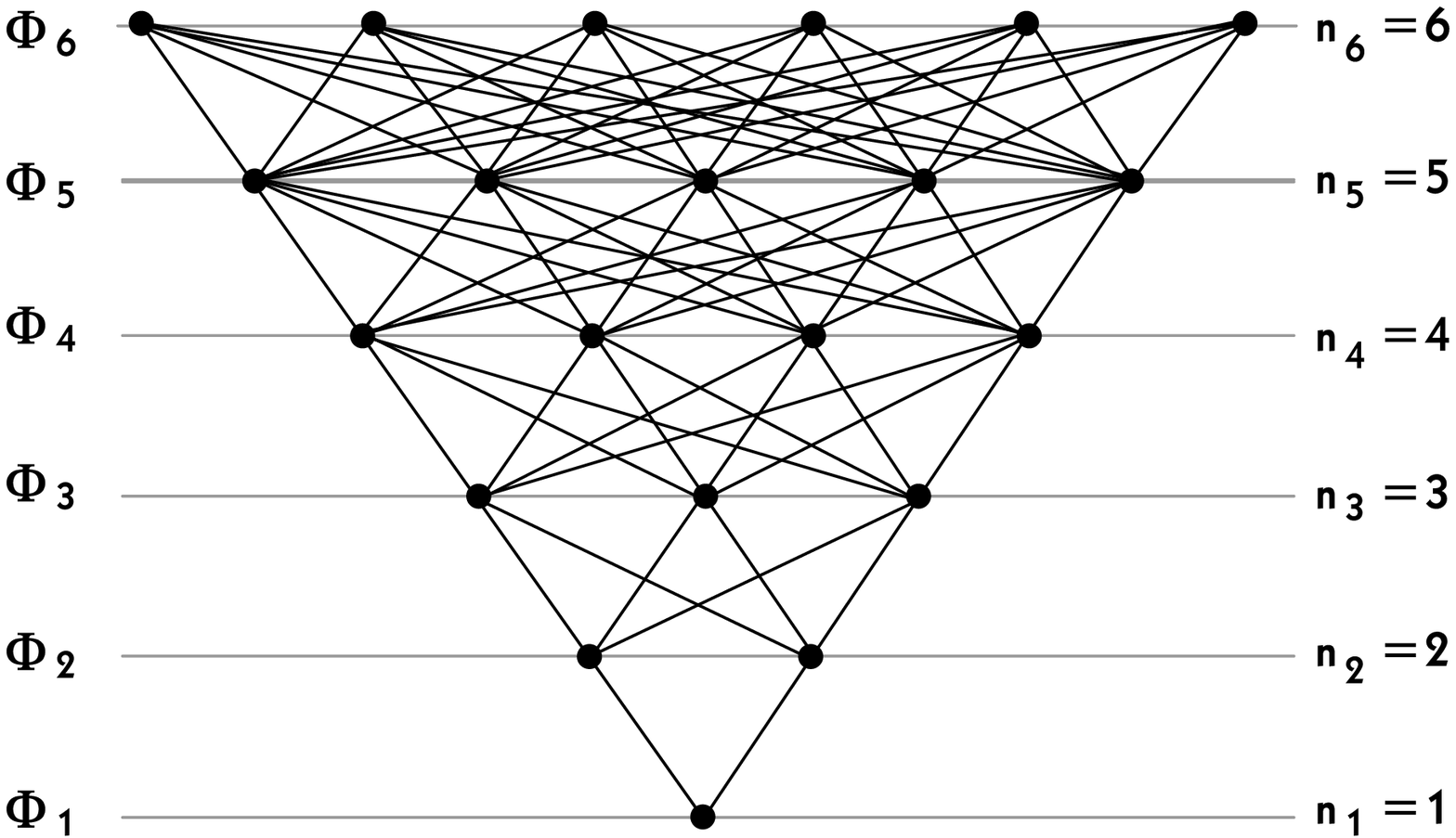}
	\caption{Display of Natural numbers' finite prime Cobweb poset}
\end{center}
\end{figure}

\begin{defn}
\emph{\cite{2}} A natural numbers' valued sequence $F = \{n_F\}_{n\geq 0}$, $F_0 =1$ is called 
cobweb-admissible iff
\begin{displaymath}
{n \choose k}_{F}\in N_0\quad for \quad k,n\in N_0.
\end{displaymath}
\end{defn}

\noindent
$F_0 = 0$ being acceptable as $0_F! \equiv F_0! = 1$. 
We adopt then the convention to call the root $\{\emptyset\}$ the "empty root".

\vspace{0.4cm}

\noindent One of the problems posed in \cite{1,2} is the one which is the subject of our note.

\vspace{0.4cm}

\noindent \textbf{The tiling problem}

\vspace{0.2cm}

\noindent Suppose now that $F$ is a cobweb admissible sequence. Under which conditions any layer $\langle\Phi_n\!\rightarrow\!\Phi_k\rangle$ may be partitioned with help of max-disjoint blocks of established type $\sigma P_m$? Find effective characterizations and/or find an algorithm to produce these partitions.

\vspace{0.2cm}

The above Kwa\'sniewski tiling problem \cite{1,2} is first of all the problem of existence of a partition an layer $\langle\Phi_k \rightarrow\Phi_n \rangle$  with max-disjoint blocks of the form $\sigma P_m$  defined as follows:

\begin{displaymath}
	\sigma P_m = C_m [F, \sigma \langle F_1, F_2, \ldots, F_m \rangle ]
\end{displaymath}

\begin{figure}[ht]
\begin{center}
	\includegraphics[width=100mm]{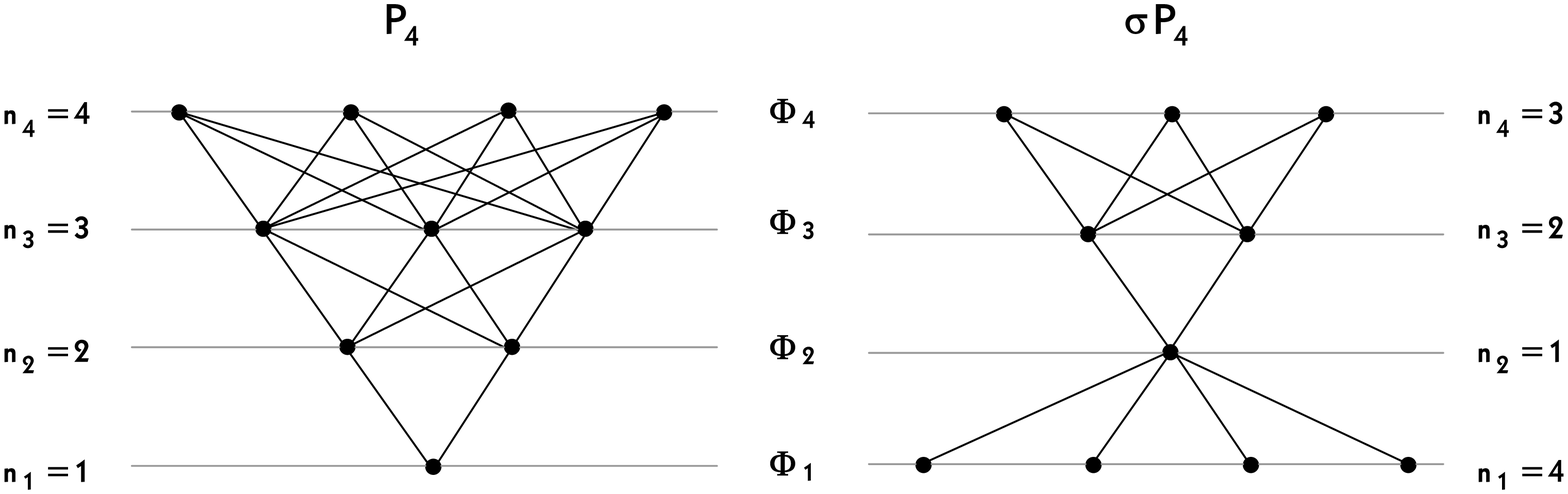}
	\caption{Display of block $\sigma P_m$ obtained from $P_m$ and permutation $\sigma$ \label{fig:permut}}
\end{center}
\end{figure}

It means that partition may contain only primary cobweb sub-posets or those obtained from primary cobweb poset $P_m$ via permuting its levels as illustrated below (Fig. \ref{fig:permut}).

\section{Example of a cobweb poset recurrent tiling algorithm - 1 (cprta1)}
 
Now we present an algorithm to create partition of any layer $\langle\Phi_k \rightarrow\Phi_n \rangle$, $k\leq\nolinebreak n$, $k,n \in \mathbb{N}\cup \{0\}$   of finite cobweb sub-poset specified by such  $F$-sequences as Natural numbers and Fibonacci numbers. We shall use the abbreviation: (cprta1) algorithm. In the following  Theorem 1  and Theorem 2 are existence theorems.

\begin{theoremn}[Natural numbers] 
Consider any layer $\langle \Phi_{k+1} \rightarrow \Phi_n \rangle$ with $m$ levels where $m=n-k$, $k\leq{n}$ and  $k,n\in\mathbb{N} \cup \{0\}$ in a finite cobweb sub-poset, defined by the sequence of \textbf{natural numbers} i.e. $F \equiv \{ n_F \}_{n\geq 0},\ n_F = n,\ \linebreak n\in\mathbb{N} \cup \{0\}$. Then there exists at least one way to partition this layer with help of max-disjoint blocks of the form $\sigma P_m$.
\end{theoremn}

\noindent Max-disjoint means that the two blocks have no maximal chain in common \cite{1,2}.

\noindent Before proving let us notice that for any $m, k \in \mathbb{N}$ such that $m+k=n$:

\begin{equation}
	n_F = m_F + k_F	\label{eq:1}
\end{equation}

\noindent where $1_F = 1$.

\begin{figure}[ht]
\begin{center}
	\includegraphics[width=60mm]{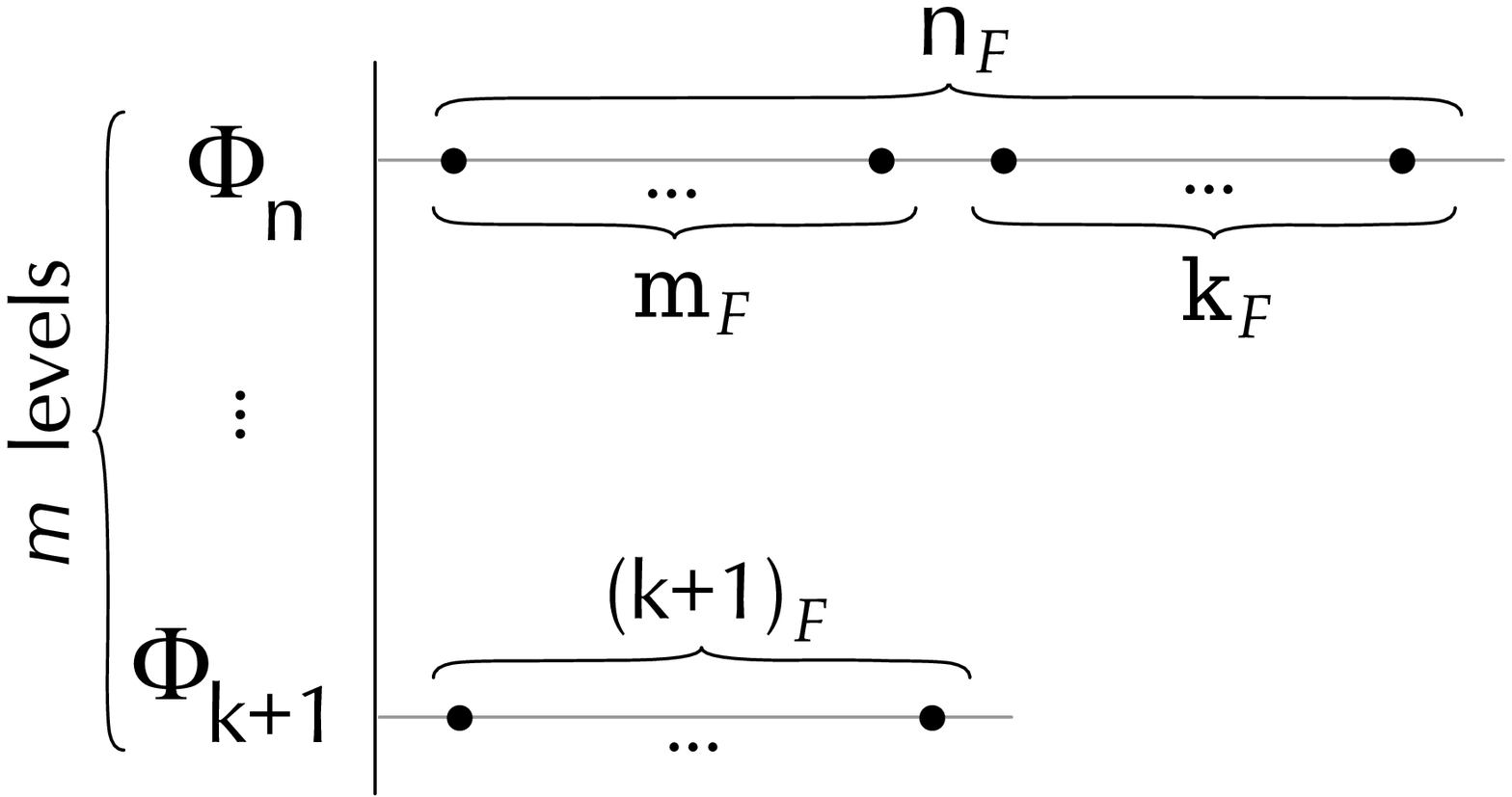}
	\caption{Picture of $m$ levels of Cobweb poset' Hasse diagram \label{fig:steep1} }
\end{center}
\end{figure}

\noindent \textbf{P}\textbf{\footnotesize{ROOF}} \textbf{(cprta1) algorithm}

\vspace{0.2cm}

\noindent \textbf{Steep 1}. There are $n_F = m_F + k_F$ vertices on the $\Phi_n$ level. Let us separate them cutting into two disjoint subsets as illustrated by the Fig.\ref{fig:steep1} and cope at first with $m_F$ vertices (Steep~2).  Then we shall cope with those $k_F$ vertices left (Steep~3).

\begin{figure}[ht]
\begin{center}
	\includegraphics[width=50mm]{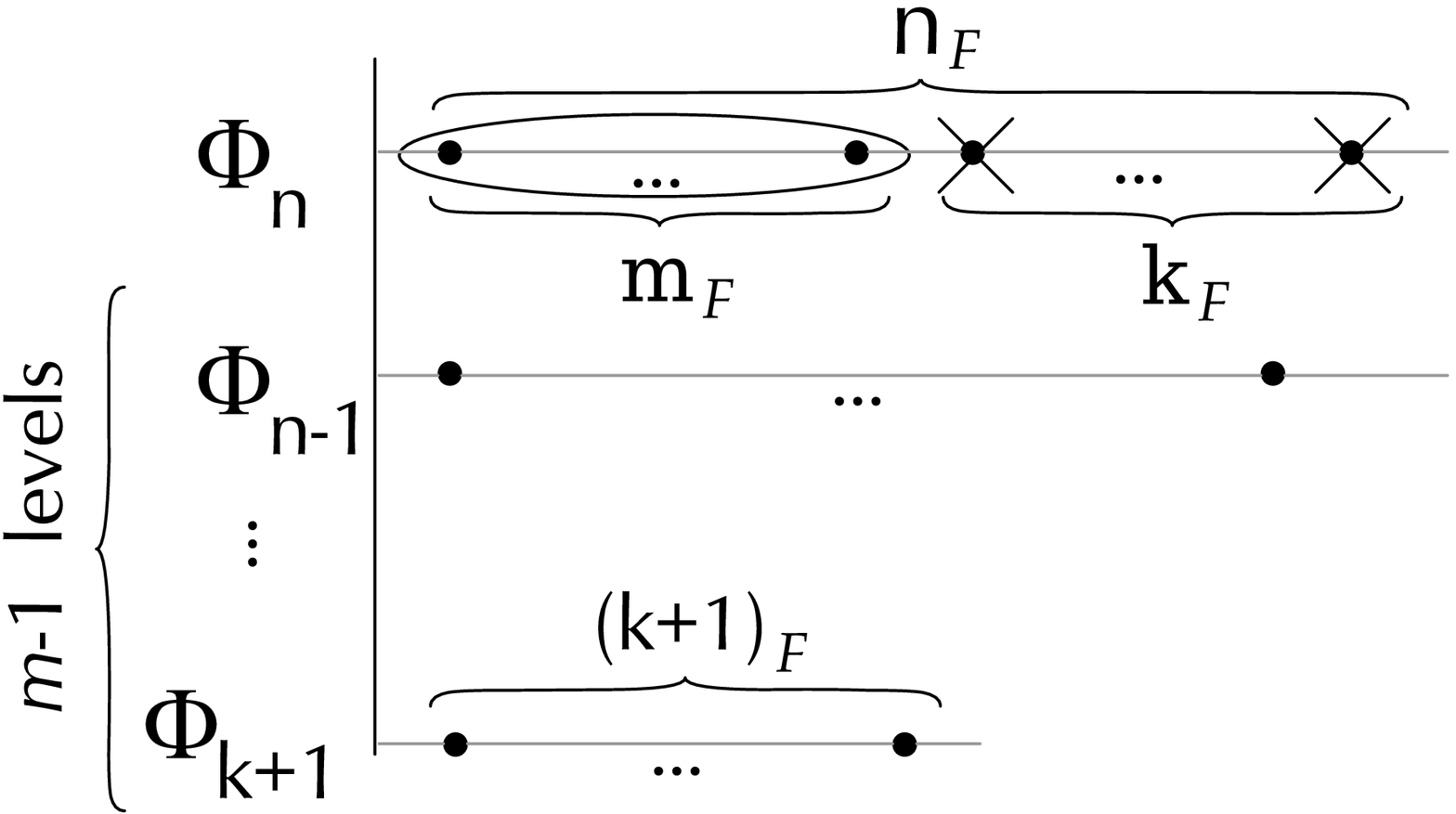}
	\caption{Picture of Steep 2 \label{fig:steep2} }
\end{center}
\end{figure}

\noindent \textbf{Steep 2}. Temporarily we have $m_F$ fixed vertices on $\Phi_n$ level to consider. Let us cover them by $m$-th level of block $P_m$, which has exactly $m_F$ vertices-leafs. What was left is the layer $\langle \Phi_{k+1} \rightarrow \Phi_{n-1} \rangle$ and we might eventually partition it with smaller max-disjoint blocks $\sigma P_{m-1}$, but we need not to do that.  See the next step.

\vspace{0.2cm}
\noindent \textbf{Steep 3}. Consider now the second complementary situation, where we have $k_F$ vertices on $\Phi_n$ level being fixed.  Observe that if we \emph{move} this level lower than $\Phi_{k+1}$ level, we obtain exactly $\langle \Phi_{k} \rightarrow \Phi_{n-1} \rangle$ layer to be partitioned with max-disjoint blocks of the form $\sigma P_m$.  This "\emph{move}" operation is just permutation of levels' order.

\vspace{0.2cm}

The layer $\langle\Phi_{k+1}\!\!\rightarrow\!\Phi_{n}\rangle$ may be partitioned with $\sigma P_m$ blocks if $\langle\Phi_{k+1}\!\rightarrow\nolinebreak\Phi_{n-1}\rangle$ may be partitioned with $\sigma P_{m-1}$ blocks and $\langle\Phi_{k}\!\!\rightarrow \Phi_{n-1} \rangle$ by $\sigma P_m$ again. Continuing these steeps by induction, we are left to prove that $\langle\!\Phi_{k}\!\rightarrow \Phi_{k} \rangle$ may be partitioned by $\sigma P_0$  blocks and $\langle \Phi_{1}\!\!\rightarrow\!\Phi_{m} \rangle$ by $\sigma P_m$ blocks which is obvious $\blacksquare$

\begin{figure}[ht]
\begin{center}
	\includegraphics[width=85mm]{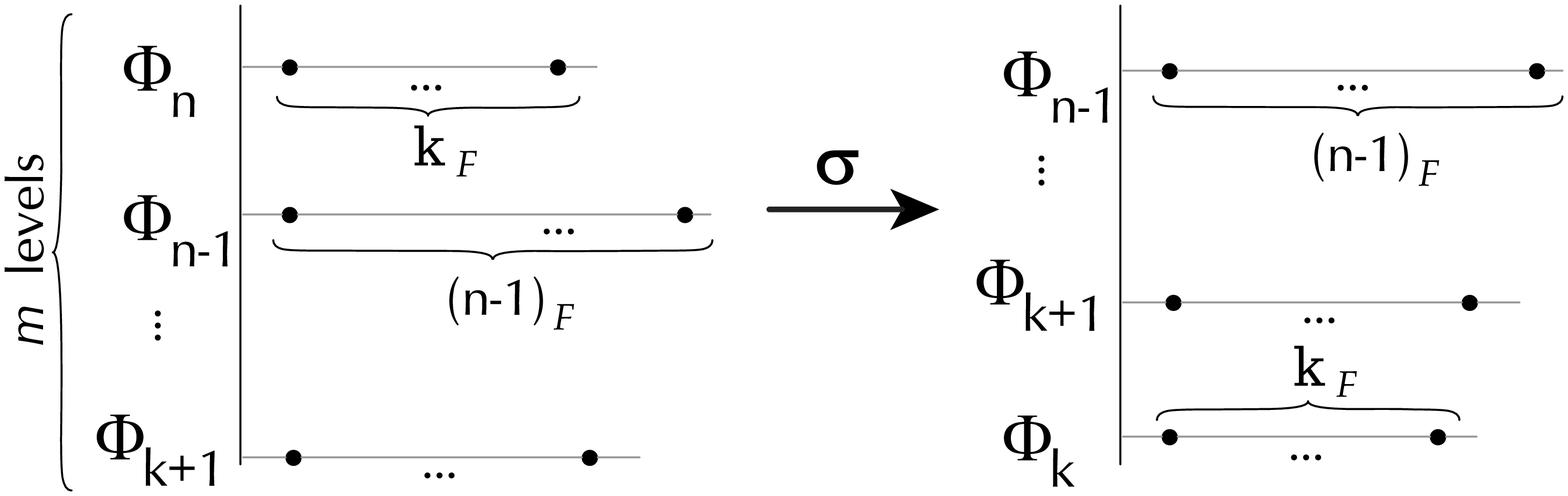}
	\caption{Picture of Steep 3 \label{fig:steep3} }
\end{center}
\end{figure}

\begin{observen}
$\\ $\emph{ We know from  \cite{1,2} (Observation 3 there) that the number of max-disjoint equip-copies of $\sigma P_m$, rooted at the same fixed vertex of $k$-th level and ending at the $n$-th level is equal to }
\end{observen}

\begin{displaymath}
	{n \choose k}_F = {n \choose m}_F
\end{displaymath}

\noindent If we cut-separate family of leafs of the layer $\langle\Phi_{k+1}\!\!\rightarrow\Phi_{n}\rangle$, as in the proof of the Theorem 1 then the number of max-disjoint equip copies of $P_{m-1}$ from the Steep 2 is equal to

\begin{displaymath}
	{n-1 \choose k}_F
\end{displaymath}

\noindent However the number of max-disjoint equip copies of $P_m$ from the Steep 3 is equal to

\begin{displaymath}
	{n-1 \choose k-1}_F
\end{displaymath}

\noindent It gives us well-known formula of Newton's symbol recurrence:

\begin{displaymath}
	{n \choose k}_F = {n-1 \choose k}_F + {n-1 \choose k-1}_F
\end{displaymath}

\vspace{0.2cm}
\noindent in accordance with  what was expected for the case  $F = \mathbb{N}$  thus illustrating  the combinatorial interpretation from  \cite{1,2}  in this particular case.  

\vspace{0.4cm}
In the next we adapt Knuth notation for "$F$-Stirling numbers" of the second kind $\Big\{ {n \atop k } \Big\}_F $ as in  \cite{2} and also in conformity with Kwa\'sniewski notation for \linebreak $F$-nomial coefficients \cite{4,1,3}.
The number of those partitions which are obtained via (cprta1) algorithm shall be denoted by the symbol $\Big\{ {n \atop k } \Big\}_F^1$.

\begin{observen} 
$\\ $\emph{ Let $F$ be a sequence matching (\ref{eq:1}). Then the number $\Big\{ {n \atop k } \Big\}_F^1 $ of different partitions of the layer $\langle\Phi_{k}\!\!\rightarrow\Phi_{n}\rangle$ where $n,k \in \mathbb{N},\ n,k \geq 1$ is equal to:}
\end{observen}

\begin{tabular}{ p{95mm}  p{10mm}}
	\begin{displaymath} 
	\bigg\{ {n \atop k } \bigg\}_F^1  =  {n_F \choose m_F}  \cdot  \bigg\{ {n-1 \atop k } \bigg\}_F^1  	\cdot  \bigg\{ {n-1 \atop k-1 } \bigg\}_F^1
	\end{displaymath}
& 
	\vspace{0.3cm}
	\begin{center} 
		$(S_N)$
	\end{center}
\end{tabular}

\noindent where $\Big\{{n\atop n}\Big\}_F^1 = \Big\{{n\atop n}\Big\}_F = 1$, 
$\Big\{{n\atop 1}\Big\}_F^1 = \Big\{{n\atop 1}\Big\}_F = 1$,
$m = n-k+1$.

\vspace{0.4cm}

\noindent \textbf{P}\textbf{\footnotesize{ROOF}}

\vspace{0.2cm}
\noindent According to the  Steep 1 of the proof of Theorem 1 we may choose on $\Phi_n$ level $m_F$ vertices out of $n_F$ ones in ${n_F \choose m_F}$ ways.  Next recurrent steps of the proof of  Theorem 1 result in formula ($S_N$) via product rule of counting. $\blacksquare$

\vspace{0.2cm}
\noindent \textbf{Note.} $\Big\{{n\atop k}\Big\}_F^1$ is not the number of all different partitions of the layer $\nolinebreak{\langle\Phi_{k}\!\!\rightarrow\!\Phi_{n}\rangle}$ i.e. $\Big\{{n\atop k}\Big\}_F \geq \Big\{{n\atop k}\Big\}_F^1 $  as computer experiments \cite{6} show. There are much more other tilings with blocks $\sigma P_m$.

\begin{figure}[ht]
\begin{center}
	\includegraphics[width=55mm]{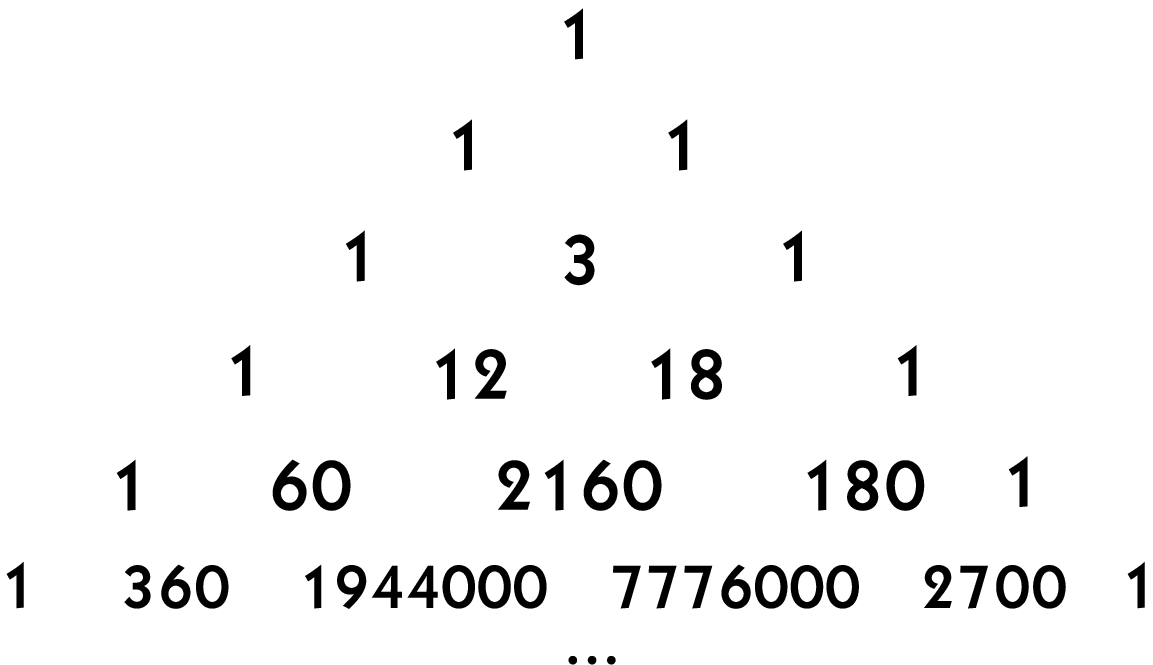}
	\caption{Natural numbers' Cobweb poset tiling triangle of $\Big\{{n\atop k}\Big\}_F^1$ \label{fig:triangle1} }
\end{center}
\end{figure}

\begin{figure}[ht]
\begin{center}
	\includegraphics[width=75mm]{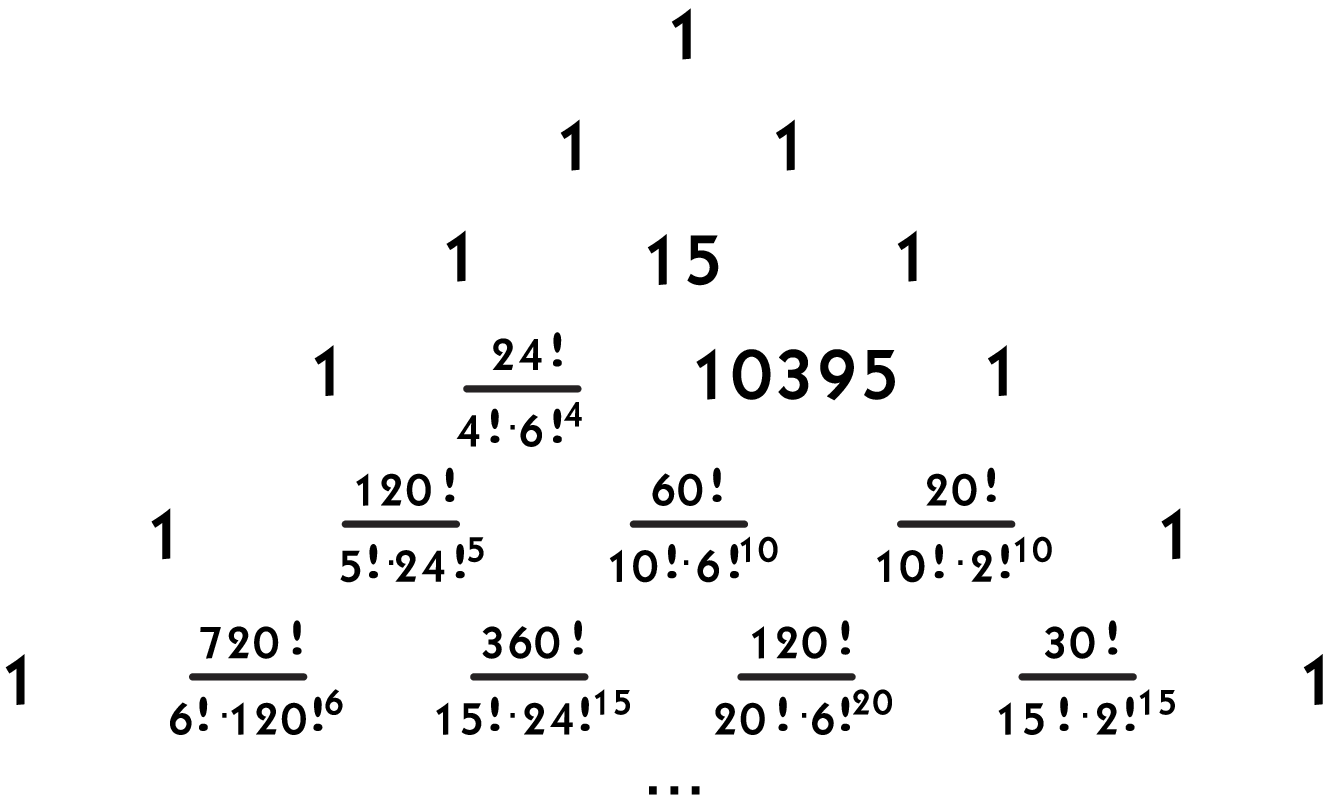}
	\caption{Kwa\'sniewski Natural numbers' cobweb poset tiling triangle of $\Big\{{\eta \atop \kappa}\Big\}_{\lambda}$ \label{fig:akk_triangle1} }
\end{center}
\end{figure}

This is to be  compared  with  Kwa\'sniewski cobweb   triangle   \cite{2} (Fig. \ref{fig:akk_triangle1}) for the infinite triangle matrix  elements

$$\Big\{{\eta \atop \kappa}\Big\}_{\lambda} = \delta_{\eta,\kappa \lambda}
\frac{\eta !}{\kappa !\lambda !^\kappa}$$
\noindent counting the number of partitions with block sizes all equal to $\lambda$.

\noindent Here $const=\lambda =  m_{F}!, m=n-k+1$ and 

$$ \eta =  n_{F}^{\underline{m}},\ \  \kappa = {n \choose k-1}_{F} $$

\noindent The inequality 
$\Big\{{n\atop k}\Big\}_F^1 \leq \Big\{{\eta \atop \kappa}\Big\}_{\lambda}$
gives us the rough upper bound for the number of tilings 
with blocks of established type  $\sigma P_m$.

\vspace{0.4cm}
\begin{theoremn}[Fibonacci numbers]
Consider any layer $\langle\Phi_{k+1}\rightarrow\Phi_n \rangle$ with $m$ levels where $m=n-k$, $k\leq n$  and $k,n\in \mathbb{N}\cup\{0\}$ in a finite cobweb sub-poset, defined by the sequence of Fibonacci numbers i.e. $F \equiv \{ n_F \}_{n\geq 0}, n_F \in \mathbb{N} \cup \{0\}$. Then there exists at least one way to partition this layer with help of max-disjoint blocks of the form $\sigma P_m$.
\end{theoremn}

The proof of the Theorem 2 for the Fibonacci sequence $F$ is similar to the proof of Theorem 1. We only need to notice that for any  $m,k \in \mathbb{N}$, $m>1$, $m+k=n$ the following identity takes place:

\begin{equation}
	n_F = (m+k)_F = (k+1)_F\cdot m_F + (m-1)_F\cdot k_F \label{eq:2}
\end{equation}
\vspace{0.4cm}
\noindent where $1_F = 2_F = 1$.

\begin{figure}[ht]
\begin{center}
	\includegraphics[width=60mm]{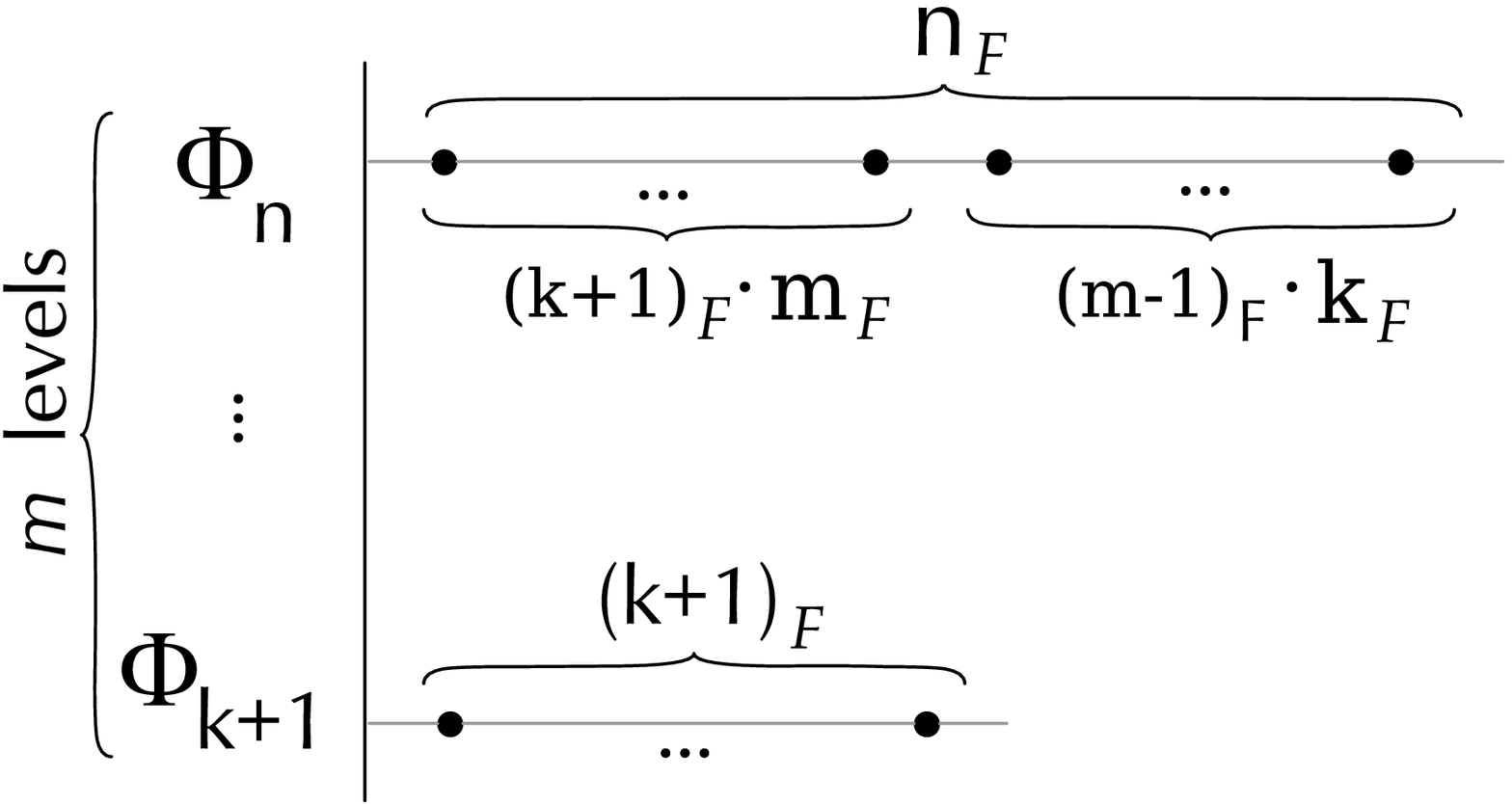}
	\caption{Picture of $m$ levels' layer of Fibonacci Cobweb graph \label{fig:steep1Fib} }
\end{center}
\end{figure}

\noindent \textbf{P}\textbf{\footnotesize{ROOF}}

\vspace{0.2cm}

\noindent The number of leafs on the $\Phi_n$ layer is the sum of two summands $\kappa\cdot m_F$  and $\mu\cdot k_F$, where $\kappa=(k+1)_F$, $\mu=(m-1)_F$, (Fig. \ref{fig:steep1Fib}) therefore as in the proof of the Theorem 1 we consider two parts.  At first we have to partition $\kappa$ layers $\langle\Phi_{k+1}\rightarrow\Phi_{n-1}\rangle$ with blocks $\sigma P_{m-1}$ and $\mu$ layers $\langle\Phi_{k}\rightarrow\Phi_{n-1}\rangle$ with $\sigma P_m$. The rest of the proof goes similar as in the case of the Theorem 1 $\blacksquare$

\vspace{0.4cm}
\noindent Theorem 2 is a generalization of Theorem 1 corresponding to $const=\kappa,\mu=1$ case.

\begin{observen}
$\\ $\emph{
The number of max-disjoint equip copies of $P_{m-1}$ which partition $\kappa$ layers $\langle\Phi_{k+1}\rightarrow\Phi_{n-1}\rangle$ is equal to}
\end{observen}

\begin{displaymath}
	\kappa{n-1 \choose k}_F = (k+1)_F{n-1 \choose k}_F
\end{displaymath}
 
However this number of max-disjoint equip copies of $P_m$ which partition $\mu$ layers $\langle\Phi_{k}\rightarrow\Phi_{n-1}\rangle$ is equal to

\begin{displaymath}
	\mu{n-1 \choose k-1}_F = (m-1)_F{n-1 \choose k-1}_F
\end{displaymath}

Therefore the  sum corresponding to the Step 2 and to the Step 3 is the well known recurrence relation for Fibonomial coefficients \cite{5,1,2,3}

\begin{displaymath}
	{n \choose k}_F = (k+1)_F{n-1 \choose k}_F + (m-1)_F{n-1 \choose k-1}_F
\end{displaymath}

\noindent in accordance with what was expected for the case $F$ being now Fibonacci sequence
thus illustrating the combinatorial interpretation from \cite{1,2} in this particular case.

\begin{observen}
$\\ $\emph{
Let $F$ be a sequence matching (\ref{eq:2}). Then the number  $\Big\{ {n \atop k } \Big\}_F^1$ of different partitions of the layer $\langle\Phi_{k}\rightarrow\Phi_{n}\rangle$ where $n,k\in\mathbb{N}, n,k\geq 1$ is equal to:}
\end{observen}

\begin{tabular}{ p{95mm}  p{1cm}}
	\begin{displaymath} 
	\bigg\{ {n \atop {k} } \bigg\}_F^1  =  \frac{F_n!}{(F_m!)^\kappa \cdot (F_{k-1}!)^\mu}
  \cdot  \bigg\{ {n-1 \atop k } \bigg\}_F^1  	\cdot  \bigg\{ {n-1 \atop k-1 } \bigg\}_F^1
	\end{displaymath}
& 
	\vspace{0.3cm}
	\begin{center} 
		$(S_F)$
	\end{center}
\end{tabular}

\noindent where
$\Big\{{n\atop n}\Big\}_F^1 = \Big\{{n\atop n}\Big\}_F = 1$, 
$\Big\{{n\atop n-1}\Big\}_F^1 = \Big\{{n\atop n-1}\Big\}_F = 1$,
$\Big\{{n\atop 1}\Big\}_F^1 = \Big\{{n\atop 1}\Big\}_F = 1$,
$\kappa=k_F, \mu=(m-1)_F, m=n-k+1$, 
$F_n! = 1\cdot 2\cdot \ldots \cdot (n_F-1)\cdot n_F$.

\vspace{0.4cm}
\noindent \textbf{P}\textbf{\footnotesize{ROOF}}

\vspace{0.2cm}
\noindent According to the Steep 1 of the proof of Theorem 2 we may choose on  $n$-th level $m_F$ vertices $\kappa$  times and next $(k-1)_F$ vertices $\mu$ times out of $n_F$ ones  in $\frac{F_n!}{(F_m!)^\kappa \cdot (F_{k-1}!)^\mu}$ ways.  Next recurrent steps of the proof of Theorem 2 result in formula  ($S_F$)  via product rule of counting $\blacksquare$

\vspace{0.4cm}
\noindent Observation 4 becomes Observation 2 once we put $const=\kappa,\mu=1$. 

\vspace{0.4cm}
\noindent {\bf Easy example} 

\begin{figure}[ht]
\begin{center}
	\includegraphics[width=70mm]{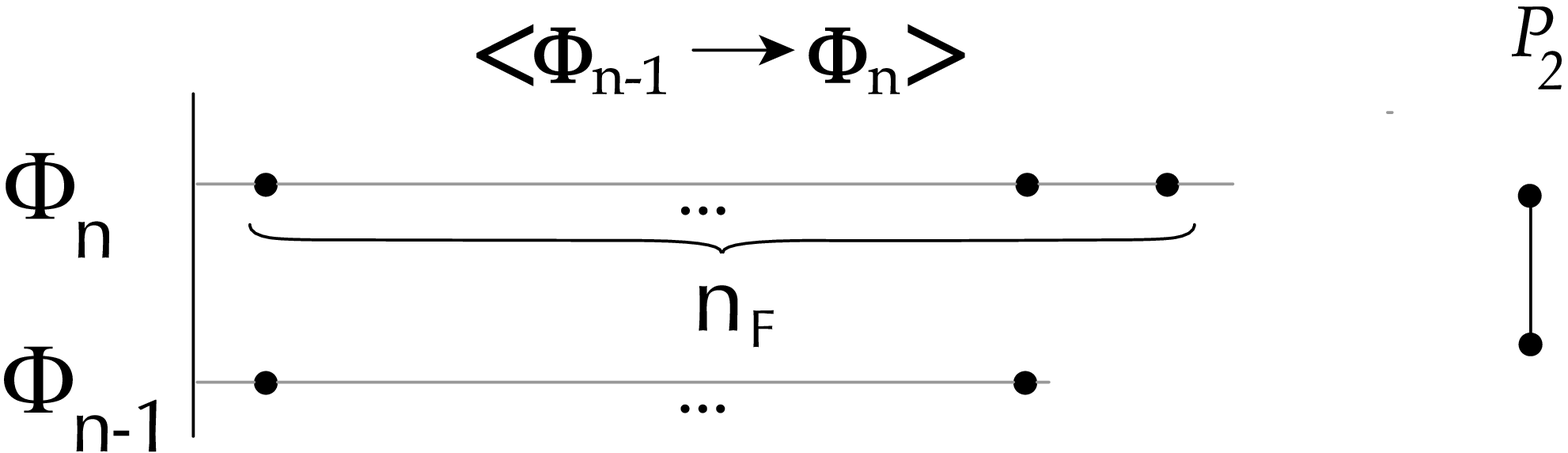}
	\caption{Easy example picture \label{fig:note11}}
\end{center}
\end{figure}

\noindent For cobweb-admissible sequences $F$ such that $1_F = 2_F = 1$,
$\Big\{\!{n\atop n-1}\!\Big\}_F^1\!\!\!=\!\!\Big\{\!{n\atop n-1}\!\Big\}_F\!\!\!=\!1$ as obviously we deal with the perfect matching of the bipartite graph which is very exceptional case (Fig. \ref{fig:note11}).

\noindent \textbf{Note.} As in the case of Natural numbers for $F$-Fibonacci numbers $\Big\{{n\atop 1}\Big\}_F^1$ is not the number of all different partitions of the layer $\langle\Phi_{k}\rightarrow\Phi_{n}\rangle$ i.e. $\Big\{{n\atop k}\Big\}_F \geq \Big\{{n\atop k}\Big\}_F^1$ as computer experiments \cite{6} show. There are much more other tilings with blocks $\sigma P_m$.

\begin{figure}[ht]
\begin{center}
	\includegraphics[width=75mm]{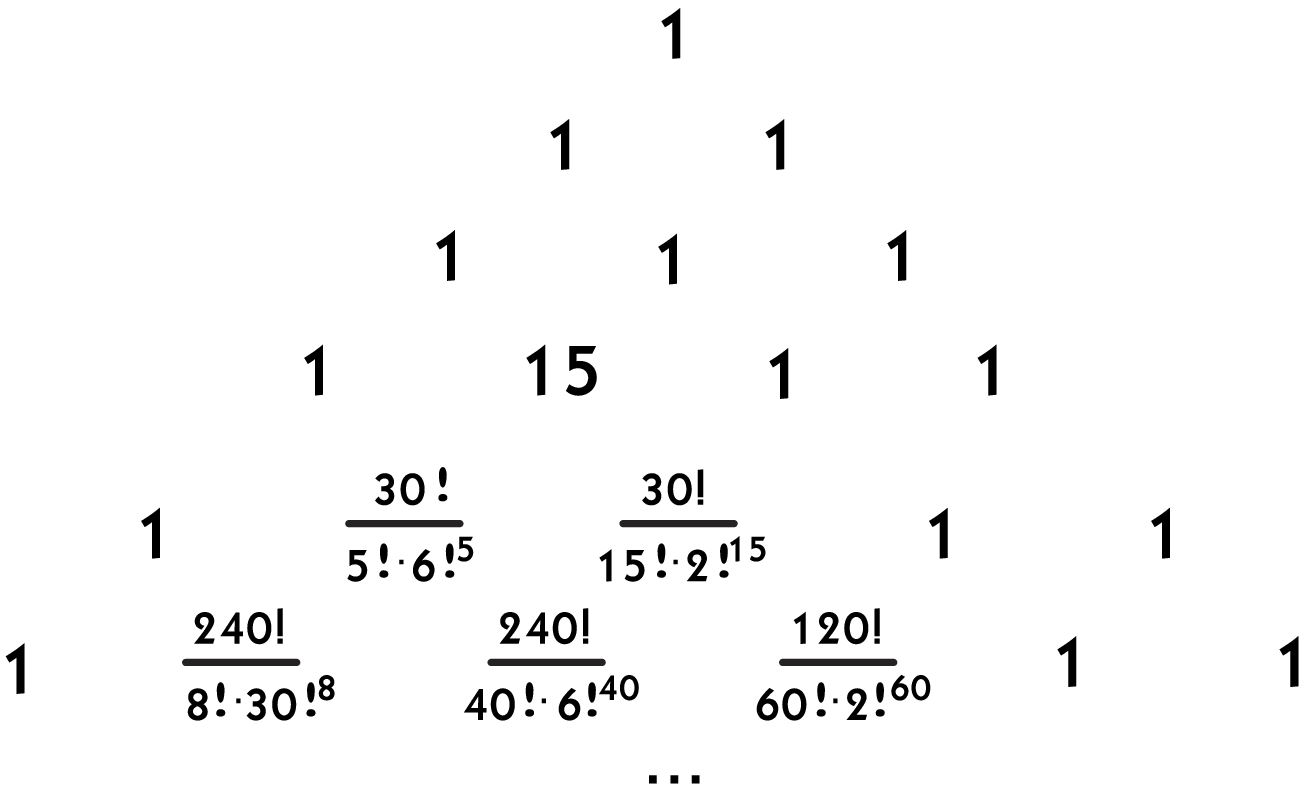}
	\caption{Kwa\'sniewski Fibonacci numbers' cobweb tiling triangle of $\Big\{{\eta \atop \kappa}\Big\}_{\lambda}$ \label{fig:akk_triangle2} }
\end{center}
\end{figure}

This is to be compared with Kwa\'sniewski cobweb triangle \cite{2} for the infinite triangle matrix elements (Fig. \ref{fig:akk_triangle2})

\begin{figure}[ht]
\begin{center}
	\includegraphics[width=60mm]{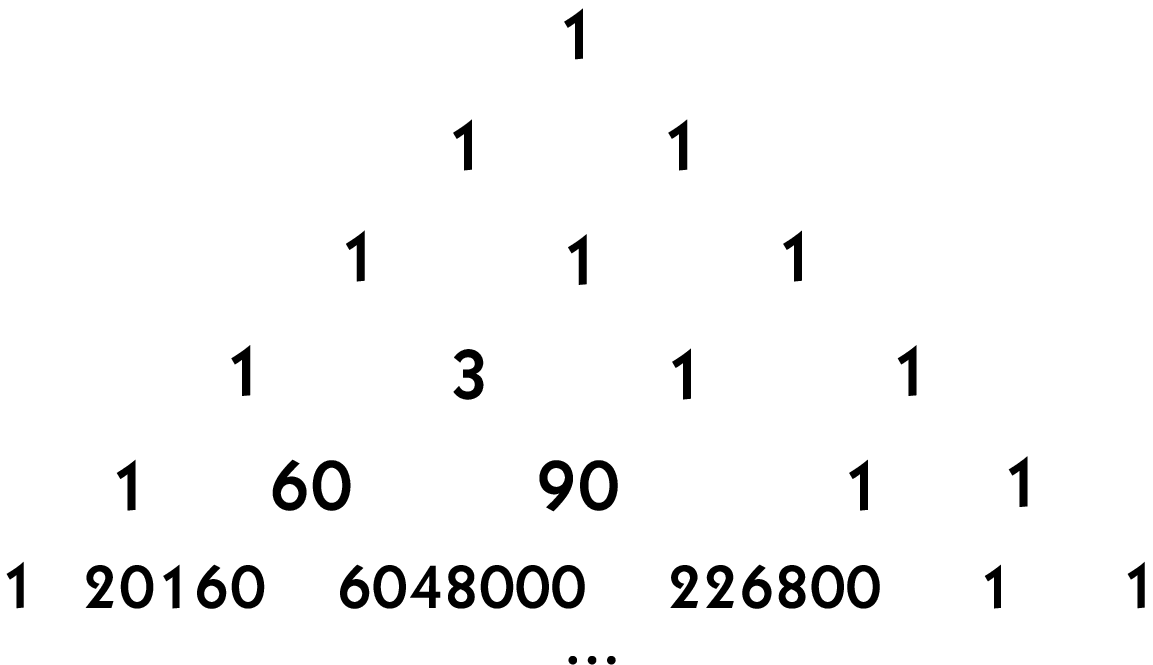}
	\caption{Fibonacci numbers' cobweb poset tiling triangle of $\Big\{{n\atop k}\Big\}_F^1$ \label{fig:triangle1_fib} }
\end{center}
\end{figure}

\section {Other tiling sequences}

\begin{defn}
The cobweb admissible sequences that designate cobweb posets with tiling are called cobweb tiling sequences.
\end{defn}

\subsection{Easy examples}
The above method applied to prove tiling existence for Natural and Fibonacci numbers relies on the assumptions (\ref{eq:1}) or (\ref{eq:2}). Obviously these are not the only sequences that do satisfy recurrences (\ref{eq:1}) or (\ref{eq:2}). There exist also other cobweb tiling sequences beyond the above ones with different initial values.
 
There exist also cobweb admissible sequences determining cobweb poset with \emph{no} tiling of the type considered in this note.

\begin{example}
$n_F = (m+k)_F = m_F + k_F$, $n\geq 1$ \emph{( $0_F$ = corresponds to one "\emph{empty root}" $\{\emptyset\}$ - compare with Definition 1 )}
\end{example} 

This might be considered a sample example illustrating the method. For example if we choose $1_F = c \in \mathbb{N}$, we obtain the class of sequences $n_F=c\cdot n$ for $n\geq 1$. Naturally layers of such cobweb posets designated by the sequence satisfying (\ref{eq:1}) for $n\geq 1$  may also be partitioned according to (cprta1).
\vspace{0.4cm}

\noindent \textbf{Example 1.5} $1_F = 1, n_F = c\cdot n, n>1$ ($0_F$ = corresponds to one "\emph{empty root}" $\{\emptyset\}$ )
 This might be considered another sample example now illustrating the "\emph{shifted}" method named (cpta2). For example if we choose $2_F = c \in \mathbb{N}$, while $1_F=1$, we obtain the class of sequences $1_F = 1$ and $n_F = c\cdot n $ for $n>1$. Layers of such cobweb posets designated by these sequences may also be partitioned.

\begin{observen} \emph{\textbf{ Algorithm (cpta2) }}

\noindent Given any (including cobweb-admissible) sequence $A\equiv\{n_A\}_{n\geq 0}$, $s\in\mathbb{N}\cup\{0\}$   let us define \emph{shift} unary operation $\oplus_s$ as follows:
\begin{displaymath}
	\oplus_s A = B, \qquad n_B = \left\{ 
	\begin{array}{ll} 
		1 				& n<s \\ 
		(n-s)_A 		& n\geq s   
	\end{array} \right. 
\end{displaymath}

\noindent where $B\equiv\{n_B\}_{n\geq 0}$. Naturally $\oplus_0$ = identity.  Then the following is true.
If a sequence $A$ is cobweb-tiling sequence then $B$ is also cobweb-tiling sequence.
\end{observen}

\vspace{0.4cm}
\noindent For example this is the case for $A=1,2,3,4,\ldots$, $\oplus_{3}A = 1,1,1,1,2,3,4,\ldots$.

\begin{example}
$n_F = m_F \cdot k_F$
\end{example} 
If we choose $1_F = c\in \mathbb{N}$, we obtain the class of sequences $n_F = c^n, n\geq0$. We can also consider more general case $n_F = \alpha\cdot m_F \cdot k_F$, where $\alpha\in\mathbb{N}$ which gives us the next class of tiling sequences $n_F = \alpha^{n-1}\cdot c^n, n\geq1, 0_F=1$ and layers of such cobweb posets can be partitioned by (cprta1) algorithm. For example: $1_F = 1, \alpha=2 \rightarrow F = 1,1,2,4,8,16,32,\ldots$ or $1_F = 2 \rightarrow F = 1,2,4\alpha,8\alpha^2,16\alpha^3,\ldots$

\newpage
\begin{example}
$n_F = (m+k)_F = (k+1)_F\cdot m_F + (m-1)_F\cdot k_F$
\end{example}
Here also we have infinite number of cobweb tiling sequences depending on the initial values chosen for the recurrence $(k\!+\!2)_F\!=\!2_F(k\!+\!1)_F\!+\!k_F, k\!\!\geq\!\!0$. For example: $1_F = 1$ and $2_F = 2 \rightarrow$ $F=1,2,5,12,29,70,169,408,985,\ldots$ 
Note that this is not shifted Fibonacci sequence as we use recurrence (2)  which depends on initial conditions adopted.
Next $1_F = 1$ and $2_F = 3 \rightarrow F=1,3,10,33,109,360,1189,\ldots$ 
Note that this is not remarcable Lucas sequence [7] (reference \cite{7} was indicated to me by A.K.Kwa\'sniewski).

\vspace{0.4cm}
Neither of sequences : shifted Fibonacci nor Lucas sequence  satisfy (2)
neither these are cobweb admissible sequences as is the case of Catalan, Motzkin, Bell or Euler numbers.

\vspace{0.4cm}
The proof of tiling existence leads to many easy known formulas for sequences, where we use multiplications of terms $m_F$ and/or $k_F$, like
$n_F = \alpha\cdot k_F$,
$n_F = \alpha\cdot m_F k_F$,
$n_F = \alpha\cdot (m\pm\beta)_F k_F$,
where $\alpha, \beta \in \mathbb{N}$, $n=m+k$ and so on.

This are due to the fact that in the course of partition's existence proving with (cprta1) partition of layer $\langle\Phi_{k+1}\!\rightarrow\!\Phi_n\rangle$ existence relies on partition's existence of smaller layers $\langle\Phi_{k+1}\!\rightarrow\!\Phi_{n-1}\rangle$ and/or $\langle\Phi_{k}\!\rightarrow\!\Phi_{n-1}\rangle$.

\vspace{0.4cm}

In what follows we shall use an at the point product of two cobweb-admissible sequences giving as a result a new cobweb admissible sequence - cobweb tiling sequences included to which the above described treatment (cprta1) applies.

\subsection{Beginnings of the cobweb-admissible sequences$\\ $ production}

\begin{defn}
Given any two cobweb-admissible sequences $A\equiv\{n_A\}_{n\geq 0}$ and $B\equiv\{n_B\}_{n\geq 0}$, their at the point product $C$ is given by
\begin{displaymath}
	A\cdot B = C \qquad C\equiv\{n_C\}_{n\geq 0},\ n_C = n_A \cdot n_B
\end{displaymath}
\end{defn}

\noindent It is obvious  that $A \cdot B = C$ is also cobweb admissible and 
\begin{displaymath}
	{n \choose k}_{A\cdot B} = \frac{n_{A}^{\underline{k}}}{k_A!} \cdot \frac{n_{B}^{\underline{k}}}{k_B!} = {n \choose k}_A \cdot {n \choose k}_B \in \mathbb{N}\cup\{0\}
\end{displaymath}

\begin{example}
Almost constant sequences $C_t$ \label{ex:const}
\begin{displaymath}
	C_t = \{n_C\}_{n\geq 0}\qquad \mathrm{where}\  const = n_C = t \in \mathbb{N}\ \mathrm{for}\ n>0, 0_F = 1.
\end{displaymath}
\end{example}

\noindent as for example $C_5 = 1,5,5,5,5,\ldots$ are trivially cobweb-admissible and cobweb tiling sequences - see next example. 

\vspace{0.4cm}
\noindent In the following $I$ denotes unit sequence $I\equiv\{1\}_{n\geq 0}$; $I\cdot{A}=A$.

\begin{example}
Not diminishing sequence $A_{c,M} $
\end{example}

\noindent If we multiply $i$-th term (where $i\geq M \geq 1, M \in \mathbb{N}$) of sequence $I$ by any constant $c \in \mathbb{N}$, then the product cobweb admissible sequence is $A_{c,M}$.
\begin{displaymath}
	A_{c,M} \equiv \{n_A\}_{n\geq 0}\qquad \mathrm{where}\ n_A = 
	\left\{ 
	\begin{array}{lr} 
		1 		& 1\leq n < M \\ 
		c 		& n\geq M   
	\end{array} \right. 
\end{displaymath}

\noindent as for example $\nolinebreak{A_{5,10} = 1,\underbrace{1,\ldots,1}_{10},5,5,5,\ldots}$ or more general example  

\noindent $\nolinebreak{A_{3,2,10} = 1,\underbrace{3,\ldots,3}_{10},6,6,6,\ldots}$ Clearly sequences of this type are cobweb admissible and  cobweb tiling sequences.

\vspace{0.4cm}
Indeed. Each of level of layer $\langle\Phi_k\!\rightarrow\!\Phi_n\rangle$ has the same or more vertices than each of levels of the block $\sigma P_m$. If not the same then the number of vertices from the block $\sigma P_m$  divides the number of vertices at corresponding layer's level. This is how   (cprta2) applies.

\vspace{0.4cm}
\noindent \textbf{Note.}
The  sequence $A_{3,2,10}$ is a product of two sequences from Example \ref{ex:const}, 
$A = 1,3,3,3,3,3,3,\ldots$ and  $B'=\oplus_{10}B= 1,\ldots,1,2,2,2,\ldots$ where $\\ B = 1,2,2,2,2,2,2,\ldots$, then $A\cdot B' = A_{3,2,10} = 1,\underbrace{3,\ldots,3}_{10},6,6,6,\ldots $ 

\begin{example}
Periodic sequence $B_{c,M}$
\end{example}

\noindent A more general example is supplied by 
\begin{displaymath}
	B_{c,M} \equiv \{n_B\}_{n\geq 0}\qquad \mathrm{where}\ n_B = 
	\left\{ 
	\begin{array}{ll} 
		1 		& M\nmid{n} \vee n=0 \\ 
		c 		& M | n  
	\end{array} \right. 
\end{displaymath}

\noindent where $c,M \in \mathbb{N}$. Sequences of above form are cobweb tiling, as for example  $B_{2,3} = \underbrace{1,1,2}_{3},1,1,2,\ldots$, $B_{7,4} = \underbrace{1,1,1,7}_{4},1,1,1,7,\ldots$ Indeed.

\vspace{0.4cm}
\noindent \textbf{P}\textbf{\footnotesize{ROOF}}
\vspace{0.2cm}
\noindent Consider any layer $\langle\Phi_k\!\rightarrow\!\Phi_n\rangle$, $k\leq n$, $k,n\in\mathbb{N}\cup\{0\}$, with $m$ levels:

\noindent For $m < M$, the block $P_m$ has one vertex on each of levels. The tiling is trivial.

\noindent For $m \geq K$, the sequence $B_{c,M}$ has a period equal to $M$, therefore any layer of $m$ levels has the same or larger number of levels with $c$ vertices than the block $\sigma P_m$, if layer's level has more vertices than corresponding level of block $\sigma P_m$ then the quotient of this numbers is a natural number i.e. $1|c$, thus the layer can be partitioned by one block $P_m$ or by $c$ blocks $\sigma P_m$ $\blacksquare$

\begin{observen}
$\\ $\emph{
The at the point product of the above sequences gives us occasionally a method to produce Natural numbers as well as expectedly other cobweb-admissible sequences with help of the following algorithm.}
\end{observen}

\vspace{0.4cm}
\noindent \textbf{Algorithm for natural numbers' generation (cta3)}

\vspace{0.4cm}
\noindent $N(s)$ denotes a sequence which first $s$ members is next Natural numbers i.e. 
$N(s)\equiv\{n_N\}_{n\geq 0}$, where $n_N = n$, for $n=1,2,\ldots,s$, $\ p, p_n$ - prime numbers.

\begin{enumerate}
	\item $N(1) = \mathbf{I} = 1,1,1,\ldots $
	\item $N(2) = N(1)\cdot B_{2,2} = 1,2,1,2,1,2,\ldots$
	\item $N(3) = N(2)\cdot B_{3,3} = 1,2,3,2,1,6,\ldots$
	\item[n.] $N(n) = N(n-1) \cdot \mathbf{X} $
\end{enumerate}

\noindent Consider $n$:
\begin{enumerate}
	\item let $n$ be prime, then $\neg \exists_{1\neq i \in [n-1]} i|n \Rightarrow n_N = 1 \Rightarrow \mathbf{X} = B_{n,n} $

	\item let $n = p^m$, $1<m \in \mathbb{N}$, then $n_N = p^{m-1} \Rightarrow \mathbf{X} = B_{p,n}$

	\item let $n = \prod_{s=1}^{u} p_{s}^{m_s}$, where $p_i \neq p_j$ for $i\neq j$, $m_i \geq 1$, $\\ i=1,2,\ldots,u$, $u>1$
	
	$\forall_{i\in[u]} p_i^{m_i} < n  \Rightarrow n_N = \mathrm{LCD}\left(\{ p_i^{m_i}:  i=1,2,\ldots,u \}\right) $
	$\\ \wedge \forall_{i\neq j} \mathrm{GCD}( p_i^{m_i},  p_j^{m_j}) = 1 $ 
	$\Rightarrow$ $n_N = \prod_{s=1}^{u} p_{s}^{m_s}$
	$\Rightarrow \mathbf{X} = \mathbf{I}$

\end{enumerate}

\noindent where lowest common denominator or least common denominator (LCD)  and  greatest common divisor (GCD)  abbreviations were used.

\vspace{0.4cm}
\noindent \textbf{Concluding}
\begin{displaymath}
	N(n) = N(n-1)\cdot B_{h_n, n} 
	\begin{array}{c}
		\scriptstyle{n \rightarrow \infty} \\
		\longrightarrow \\
		\ 
	\end{array}
	\mathbb{N}
\end{displaymath}
\begin{displaymath}	
	h_n = \left\{ 
	\begin{array}{ll}
	p & \quad n = p^m,\quad \mathbb{N} \ni m \geq 1  \\
	1 & \quad n = \prod_{s=1}^{u>1} p_{s}^{m_s},\quad \mathbb{N} \ni m_s \geq 1
	\end{array}
	\right.
\end{displaymath}

\vspace{0.2cm} \noindent 
while  $\{h_n\}_{n\geq 1} = 1,2,3,2,5,1,7,2,3,1,11,1,13,1,1,2,17,\ldots$

\begin{figure}[ht]
\begin{center}
	\includegraphics[width=100mm]{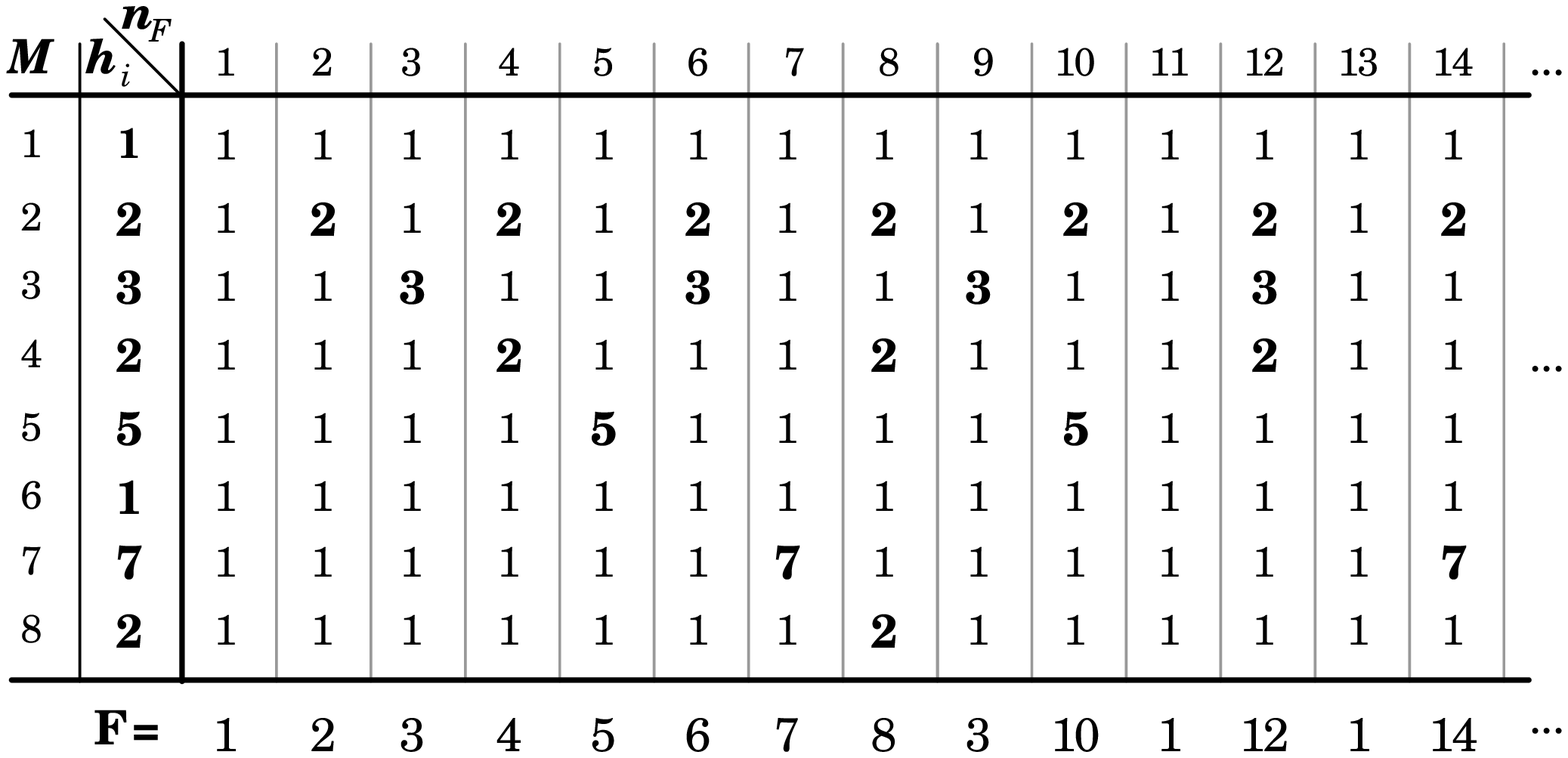}
	\caption{Display of eight steeps of algorithm (cta3) \label{fig:gennat}}
\end{center}
\end{figure}

\vspace{0.2cm}
As for the Fibonacci sequence we expect the same statement to be true for $n\rightarrow\infty$ bearing in mind those properties of  Fibonacci numbers which make them an effective tool in Zeckendorf representation of natural numbers. For the Fibonacci numbers the would be sequence $\{h_n\}_{n\geq 1}$  is given by $\{h_n\}_{n\geq 1}=1,1,2,3,5,4,13,7,17,11,89,6,\ldots$

\vspace{0.2cm}
We end up with general observation - rather obvious but important to be noted.

\newpage
\begin{theoremn}
Not all cobweb-admissible sequences are cobweb tiling sequences. 
\end{theoremn}

\noindent \textbf{P}\textbf{\footnotesize{ROOF}}

\begin{figure}[ht]
\begin{center}
	\includegraphics[width=110mm]{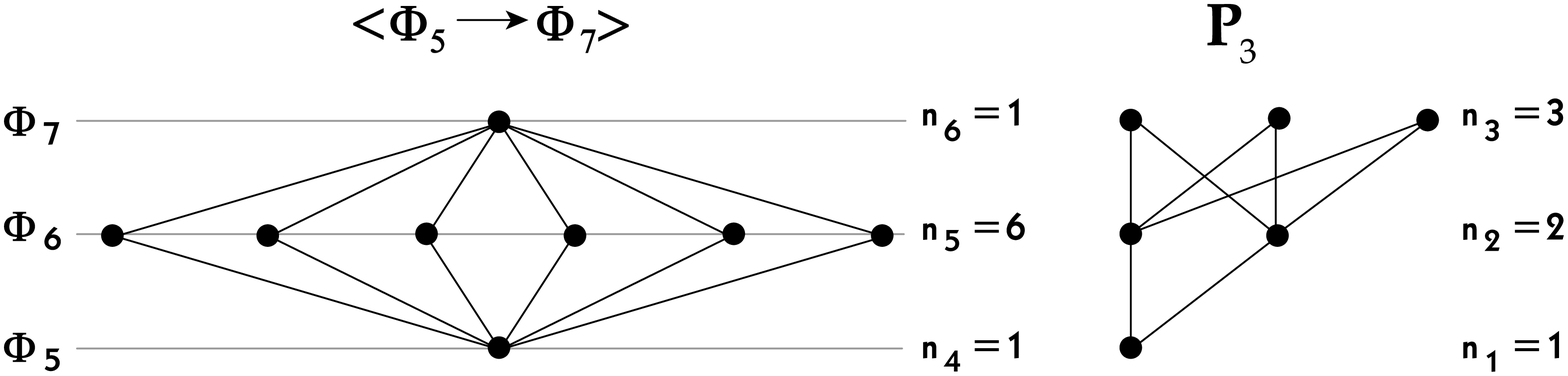}
	\caption{Picture proof of Theorem 3 \label{fig:contr} }
\end{center}
\end{figure}

\noindent It is enough to give an appropriate example. Consider then a cobweb-admissible sequence $F=A\cdot B=1,2,3,2,1,6,1,2,3,\ldots$, where $A=1,2,1,2,1,2\ldots$  and $B=1,1,3,1,1,3,\ldots$ are both cobweb admissible and cobweb tiling.
Then the layer $\langle\Phi_5\!\rightarrow\!\Phi_7\rangle$ can not be partitioned with blocks $\sigma P_3$ as the level $\Phi_5$ has one vertex, level $\Phi_5$ has six while $\Phi_5$ has one vertex again (Fig \ref{fig:contr}). 

\vspace{0.4cm}
\noindent \textbf{Corollary}
\emph{The at the point product of two tiling sequences does not need to be a tiling sequence.}

\vspace{0.4cm}
However  for  $A=1,2,1,2,\ldots$ and $B=1,1,3,1,1,3,\ldots$ cobweb tiling sequences their product $F=A\cdot B=1,2,3,2,1,6,1,\ldots$ is not a cobweb  tiling sequence.

 A natural question - enquire is anyhow still ahead \cite{1,2}. Find the effective characterizations and or algorithms for a cobweb admissible sequence to be a cobweb tiling sequence.

\vspace{0.4cm}
\noindent \textbf{Acknowledgements}

\vspace{0.4cm}
I would like to thank Professor A. Krzysztof Kwa\'sniewski - who initiated my interest in his cobweb poset concept - for his very helpful comments, suggestions, improvements and corrections of this note.


\end {document}